%
%
%
%
\magnification=1100
\font\bfeins=cmbx12
\input amstex
\documentstyle{amsppt}
\pagewidth{16.42truecm}
\TagsOnRight
\define\abs{\par\vskip 0.3cm\goodbreak\noindent}
\define\CC{C\!\!\!\!I}
\define\RR{I\!\!R}

\define\pair#1#2{{(#1,#2,\langle\cdot,\cdot\rangle)}}
\define\m{\roman m}
\define\M{\roman M}
\define\opop#1{{#1}^{\roman{op}}_{\roman{op}}}
\define\E{\bold 1\!\! \roman l}
\define\id{\roman{id}}
\define\lfl{\leaders\hbox to 1em{\hss \hss}\hfill}
\define\oml#1{{}_{#1}\omega}
\define\omr#1{\omega_{#1}}
\define\lact{{\triangleright}}
\define\ract{{\triangleleft}}
\define\nl{\par\noindent}
\topmatter
\title \bfeins
Pairing and Quantum Double
of Multiplier Hopf Algebras
\endtitle
\author Bernhard Drabant\\
     and\\
        Alfons Van Daele
\endauthor
\thanks B.\ D.\ is supported by a grant of the KU Leuven. Address after
January 1997: Institute for Theoretical Physics,
\phantom{iiiii}University of Valencia,
Spain.
\endthanks
\address Department of Mathematics,
Katholieke Universiteit Leuven
\nl
\hbox to .43cm{\ \hfill}Celestijnenlaan 200B,
B-3001 Heverlee, Belgium
\endaddress
\email bernhard.drabant\@wis.kuleuven.ac.be
\nl
\hbox to 2.62cm{\ \hfill}alfons.vandaele\@wis.kuleuven.ac.be
\endemail
\date June 1996
\enddate
\abstract
We define and investigate pairings of multiplier Hopf algebras.
It is shown that two dually paired regular multiplier
Hopf ($*$-)algebras $A$ and $B$ yield a quantum double
multiplier Hopf ($*$-)algebra which is again regular.
Integrals on $A$ and $B$ induce an integral on the quantum double.
The results generalize pairing and quantum double construction
from ordinary Hopf algebras to multiplier Hopf algebras.
\endabstract
\toc
\widestnumber\head{3}
\specialhead{}Introduction\endspecialhead
\head 1. Preliminaries on Multiplier Hopf Algebras
\endhead
\head 2. Pairing of Multiplier Hopf Algebras
\endhead
\head 3. The Quantum Double
\endhead
\specialhead{}Appendix\endspecialhead
\endtoc
\endtopmatter

\document
\baselineskip=14pt
\advance\hoffset by -0.5truecm
\leftheadtext\nofrills{B.\ Drabant and A.\ Van Daele}
\rightheadtext\nofrills{Pairing and Quantum Double of Multiplier Hopf
   Algebras}
%
%

\specialhead{}\centerline{Introduction}\endspecialhead
\abs
The non-commutative generalization of the abelian
$C^*$-algebra of continuous complex functions over a compact group are
the so-called compact quantum groups or compact quantum group algebras
\cite{Wo1,DK}. The notion of
a Hopf algebra enters the construction of such objects.
A multiplier Hopf algebra $A$ is a not necessarily unital generalization
of Hopf algebras where the image of the comultiplication $\Delta$
is contained in the multiplier algebra $\M(A\otimes A)$, instead of
$A\otimes A$ \cite{VD2}. If $(A,\Delta)$ has an integral \cite{VD3,Swe}
and is regular
-- i.e.\ also the co-opposite multiplier Hopf algebra
$(A,\Delta^{\roman{op}})$ exists--
then the dual $(\hat A,\hat \Delta)$ is again a regular multiplier Hopf
algebra and has an invariant integral \cite{VD3}. It is also shown in
\cite{VD3} that the dual of $(\hat A,\hat \Delta)$ is canonically
isomorphic to $(A,\Delta)$.
So, in this case, duality can now be described within the same category.
\nl
For instance the algebra $C_c(G)$ of (continuous) complex functions
with compact support on a discrete group $G$ is a multiplier Hopf algebra
in a canonical way \cite{VD3,VD4}.
Also the discrete quantum groups \cite{ER,VD4}
(as well as the compact quantum groups) are multiplier Hopf algebras.
And therefore the duality of discrete quantum groups and compact
quantum groups \cite{ER,VD4} turns out to be the duality in the category of
multiplier Hopf algebras.
\nl
There is good hope to extend the notion of (regular) multiplier Hopf
algebras (with integral) to a topological version. This could
serve as fruitful starting point for a systematic definition of
locally compact quantum groups. It seems that the new theory then
containes all the special
examples existing so far in the literature \cite{MN,PW,Wo2} including the
locally compact groups. Also there duality and the existence of a positive
integral or Haar measure will play an important r\^ole.
\nl
In the present article we are interested in the more general notion of
pairing of (regular) multiplier Hopf algebras. The dual pairing
of $A$ and $\hat A$ will be seen to be a special case. This has already
been announced in \cite{VD3}. We show that two
dually paired multiplier Hopf algebras admit the construction of
a quantum double object which is again a multiplier Hopf algebra.
Regularity and $*$-property as well as the existence of an
integral can be proven also for the quantum double.
Hence we are able to construct a quantum double multiplier Hopf
algebra within the same category. This procedure yields
further interesting, non-trivial examples of multiplier Hopf algebras.
\nl
The results of this paper generalize the well known properties of
Hopf algebra pairing \cite{Ma1,VD1} and the construction of a quantum
double out of two dually paired Hopf algebras \cite{Dri,Ma2,VD1}.
Although there
is an obvious loss of categorical symmetry in the defining equations
passing from Hopf algebras to multiplier
Hopf algebras many features of the theory of Hopf algebras can be
extended to the multiplier Hopf algebra setting.
One reason for this is the fact that the defining Hopf
relations generalize to the level of the multipliers.
However it is not yet
clear, for instance, if for the pairing $(A,\hat A)$ the quantum double
multiplier Hopf algebra $\Cal D(A)$ can be reconstructed from a
category of modules as in the usual Hopf algebra case \cite{JS,Ma3}.
This is one of the open questions which are currently under investigation.
\nl
In Chapter 1 we repeat the main definitions and results on multipliers
and multiplier Hopf algebras and provide several lemmas and propositions
which are used in the sequel. Chapter 2 introduces the notion of
(pre-)pairings of multiplier Hopf algebras. The definition of
so-called multiplier Hopf algebra pre-pairing leads to several equivalent
conditions which serve as additional axioms for the definition of
multiplier Hopf algebra pairing. The ordinary Hopf algebra pairing is a
special case
of this construction. Using the results of Chapter 2 we construct
in Chapter 3 the quantum double of a dually paired couple of
multiplier Hopf algebras $(A,B)$. We will prove that the quantum double
is again a regular multiplier Hopf algebra. There exists an integral on
the quantum double if those exist on $A$ and $B$. 
If $A$ and $B$ are multiplier Hopf $*$-algebras we prove that
the quantum double has a $*$-structure.
In many of our calculations we use a
``generalized Sweedler notation'' which will be outlined in the Appendix.
\abs\abs
\head 1. Preliminaries on Multiplier Hopf Algebras\endhead
\abs
Henceforth we work with modules over the field $k=\CC$ or
$k=\RR$. By an associative algebra $A$ (over $k$)
we mean an algebra which {\it need not} contain a unit.
Hence this notion
is more general than the one for unital algebras. We suppose that
all algebras under consideration have a non-degenerate product,
i.e.\ $ab=0$ for all $a\in A$ implies $b=0$ and from
$ab=0$ for all $b\in A$ it follows that $a=0$.
With $A$ and $B$ two non-degenerate algebras the tensor algebra
$A\otimes B$ is obviously non-degenerate, too.
\nl
A multiplier $\rho=(\rho_1,\rho_2)$ of the algebra $A$ is
a pair of linear mappings in $\roman{End}_k(A)$
such that $\rho_2(a)b=a\rho_1(b)$
for all $a,b\in A$. The set of multipliers of $A$ will be denoted by
$\M(A)$. It is a unital algebra which contains $A$ as essential ideal
through the embedding $a\hookrightarrow (a\cdot,\cdot a)$.
Hence $\rho\cdot a=(\rho_1(a)\cdot,\cdot\rho_1(a)) \equiv\rho_1(a)$
and $a\cdot\rho=(\rho_2(a)\cdot,\cdot\rho_2(a))\equiv\rho_2(a)$
for all $\rho\in\M(A)$ and $a\in A$. Therefore we will frequently
use the identification $a\cdot\rho=\rho_2(a)$ and $\rho\cdot a=\rho_1(a)$.
If $A$ is unital then $A=\M(A)$. 
If $A$ is a $*$-algebra then $\M(A)$ is a $*$-algebra through
$\rho^*=(\rho_2^*,\rho_1^*)$ where $\psi^*(a):=\psi(a^*)^*$ for any
$a\in A,\psi\in \roman{End}_k(A)$.
Since the multiplication of $A$ is supposed to be non-degenerate
a multiplier $\rho=(\rho_1,\rho_2)$ of $A$ is uniquely determined by
its first or second component. 
For a tensor product of two algebras
$A$ and $B$ one obtains the canonical algebra embeddings
$$A\otimes B\hookrightarrow\M(A)\otimes\M(B)\hookrightarrow\M(A\otimes B)
   \,.\tag 1.1
$$
We often work with extensions of algebra morphisms and module maps
without mentioning it explicitely.
In the following we will outline this notation. We
refer the reader to this exposition whenever she or he suspects to
meet extensions in the course of the paper. 
\nl
Let $A$ and $B$ be algebras, and $\varphi :A\to \M(B)$
be an algebra morphism. Then $\varphi $ is called non-degenerate
algebra morphism if $B=\roman {span}\{\varphi(a)b\mid a\in A,\ b\in B\}
=\roman {span}\{b\varphi(a)\mid a\in A,\ b\in B\}$.
Analogous conditions hold for non-degenerate $*$-algebra morphisms.
We call an $A$-left module $X$ non-degenerate with respect to $A$ if the
module map $\mu:A\otimes X\to X$ is surjective
and if $\mu(a\otimes x)=0$ for all $a\in A$ implies $x=0$. A similar
definition holds for $A$-right modules. The following propositions can
now be proved in a similar way as outlined in \cite{VD2}.

\proclaim{Proposition 1.1}
Any non-degenerate algebra morphism has a unique extension to an algebra
morphism $\varphi:\M(A)\to\M(B)$.\lfl$\square$
\endproclaim

\proclaim{Proposition 1.2} Let $A$ and $B$ be algebras,
and $B$ be an non-degenerate $A$-left module through
the module map $\mu:A\otimes B\to B$
Then there exists a unique extension $\mu:\M(A)\otimes B\to B$ rendering
$B$ an $\M(A)$-left module.\lfl$\square$
\endproclaim
\nl
These notions of non-degeneracy are automatic for unital algebras.
We will now give the definition of multiplier Hopf algebras as they were
introduced in \cite{VD2}.

\definition{Definition 1.3}
Let $A$ be an algebra. An algebra
morphism $\Delta:A\to \M(A\otimes A)$ is called a comultiplication on $A$
if for all $a,a'\in A$ 
$$\left.\aligned
T_1(a\otimes a')&:= \Delta(a)(\E\otimes a')\\
T_2(a\otimes a')&:= (a\otimes\E)\Delta(a')
\endaligned\right\}
\in A\otimes A
\tag 1.2
$$
and if the linear mappings $T_1,\, T_2:A\otimes A\to A\otimes A$ obey the
relation
$$(T_2\otimes\id)\circ(\id\otimes T_1)=
  (\id\otimes T_1)\circ(T_2\otimes\id)\,.
  \tag 1.3
$$ 
If $T_1$ and $T_2$ are bijective then the pair $(A,\Delta)$ is called a
multiplier Hopf algebra or shortly MHA. If $A$ is a $*$-algebra we
demand $\Delta $ to be a $*$-algebra homomorphism.
The multiplier Hopf algebra $(A,\Delta)$ is called regular if
in addition $A_{\roman{op}}:=(A,\Delta^{\roman{op}})$ is a multiplier Hopf
algebra, where
$\Delta^{\roman{op}}$ is the opposite comultiplication,
$\Delta^{\roman{op}}(a)(b\otimes c)=\tau(\Delta(a)(c\otimes b))$ for
$a,b,c\in A$, and
henceforth $\tau:A \otimes A \to A\otimes A$ is the usual
tensor transposition.
\enddefinition

\remark{Remark 1}
Equation (1.3) replaces the coassociativity of the comultiplication of
ordinary Hopf algebras. The definition of a multiplier Hopf algebra
however guarantees that
the comultiplication is coassociative in the sense that
$$(\id\otimes\Delta)\circ\Delta=(\Delta\otimes\id)\circ\Delta
  :\M(A)\to\M(A\otimes A\otimes A)\,.
  \tag 1.4
$$
This fact is used in particular in the appendix to define a ``generalized
Sweedler notation'' which will be helpful in the calculations
of many proofs in the paper.
\endremark
\abs
From \cite{VD2} it is known that multiplier Hopf algebras automatically
possess a unique counit $\varepsilon$ and an antipode $S$ such that
$\varepsilon(a)\,a'=\m\circ T_1^{-1}(a\otimes a')$ and
$S(a)\,a'=(\varepsilon\otimes\id)\circ T_1^{-1}(a\otimes a')$.
For regular MHA's the antipode is bijective and
$S_{\roman{op}}=S^{-1}$, $\varepsilon_{\roman{op}}=\varepsilon$.
In this case the
corresponding mappings $T_{\roman{op}\,1}$ and $T_{\roman{op}\,2}$
for $A_{\roman{op}}$
can be expressed as follows.
$$\aligned
T_{\roman{op}\,1}&=(\id\otimes S^{-1})\circ\tau\circ T_2^{-1}\circ\tau\circ
      (\id\otimes S)\\
T_{\roman{op}\,2}&=(S^{-1}\otimes\id)\circ\tau\circ T_1^{-1}\circ\tau\circ
      (S\otimes\id)\,.
\endaligned
\tag 1.5
$$      
\nl
For a multiplier Hopf algebra $(A,\Delta)$ and any linear
functional $\omega\in A'$ one can define a multiplier
$(\id\otimes\omega)\Delta(a)\in \M(A)$
for any $a\in A$ according to
$[(\id\otimes\omega)\Delta(a)]\cdot a':= 
(\id\otimes\omega)(\Delta(a)\cdot (a'\otimes \E))$ and
$a'\cdot[(\id\otimes\omega)\Delta(a)]:= 
(\id\otimes\omega)((a'\otimes \E)\cdot\Delta(a))$ for all $a'\in A$.
Analogous results hold for
$(\omega\otimes\id)\Delta(a)$.
In the same manner the following statements can be proven easily.

\proclaim{Lemma 1.4}
Let $(A,\Delta)$ be a multiplier Hopf algebra and $\omega\in A'$
be a linear functional of $A$. Then 
$(\omega\otimes\id\otimes\id)(\id\otimes\Delta)\Delta(a)$
is a multiplier in $\M(A\otimes A)$ for all $a\in A$ and it holds
$$(\omega\otimes\id\otimes\id)(\id\otimes\Delta)\Delta(a)
  =\Delta((\omega\otimes\id)\Delta(a))\,.
  \tag 1.6
$$
Analogously one obtains the multiplier identity
$$(\id\otimes\id\otimes\omega)(\Delta\otimes\id)\Delta(a)=
\Delta((\id\otimes\omega)\Delta(a))\,.
\tag 1.7
$$
\ \lfl$\square$
\endproclaim
\abs\abs
\head 2. Pairing of Multiplier Hopf Algebras\endhead
\abs
In this chapter we consider bilinear functionals between 
regular multiplier Hopf algebras. We give the definitions of
pre-pairing and pairing of multiplier Hopf algebras, and we deduce
results which are necessary for the investigation of quantum doubles
of MHAs. It will be seen that the pairing of two ordinary Hopf algebras
and the pairing of a regular multiplier Hopf algebra with non-trivial
integral with its dual $\hat A$ \cite{VD3} are special cases of MHA
pairings.

\definition{Definition 2.1}
Let $A$ and $B$ be two regular multiplier Hopf algebras,
and $\langle\cdot,\cdot\rangle:A\otimes B\to k$ be a linear mapping.
Define for all $a\in A$ and $b\in B$ the linear functionals
${\oml a:= \langle a,\cdot\rangle\in B'}$ and
$\omr b:= \langle\cdot,b\rangle\in A'$, and assume that they obey the
following properties, where $a'\in A$ and $b'\in B$.
\roster
\item $(\oml a \otimes \id)\Delta(b) \in B$
      and
      $(\id\otimes\oml a)\Delta(b)\in B$,
\item $(\omr b \otimes \id)\Delta(a) \in A$
      and
      $(\id\otimes\omr b)\Delta(a)\in A$,
\item $\oml a(\id\otimes\oml {a'})\Delta(b)=\oml {a'}(\oml a\otimes\id)
      \Delta(b)=\oml{a\cdot a'}(b)$,
\item $\omr b(\id\otimes\omr{b'})\Delta(a)= \omr{b'}(\omr b\otimes\id)
      \Delta(a)=\omr{b\cdot b'}(a)$.
\endroster
Then $\pair AB$ is called a multiplier Hopf algebra
pre-pairing. 
The pre-pairing $\pair AB$ is called non-degenerate
if $A$ and $B$ are dual with respect to the bilinear form
$\langle\cdot,\cdot\rangle$.
\enddefinition

\remark{Remark 2} The usual Hopf algebra pairing is a special case
of a multiplier Hopf algebra pre-pairing because (1) and (2) of
Definition 2.1 are trivially fulfilled, and for instance
$\oml a(\id\otimes \oml {a'})\Delta(b)={\langle a\otimes a',\Delta(b)
\rangle}
=\oml {a'}(\oml a\otimes \id)\Delta(b)=\langle a a',b\rangle=
\oml{a a'}(b)$.
\endremark
\abs
For any multiplier Hopf algebra pre-pairing $\pair AB$
we can therefore define the linear mappings
$$\align
\mu^l_{A,B}&:\cases A\otimes B&\to B\\
                    a\otimes b&\mapsto (\id\otimes\oml a)\Delta(b)
             \endcases
            \tag 2.1 \\
\mu^r_{A,B}&:\cases B\otimes A&\to B\\
                    b\otimes a&\mapsto (\oml a\otimes\id)\Delta(b)
             \endcases
            \tag 2.2 \\
\mu^l_{B,A}&:\cases B\otimes A&\to A\\
                    b\otimes a&\mapsto (\id\otimes\omr b)\Delta(a)
             \endcases
            \tag 2.3 \\            
\mu^r_{B,A}&:\cases A\otimes B&\to A\\
                    a\otimes b&\mapsto (\omr b\otimes\id)\Delta(a).
             \endcases
            \tag 2.4
\endalign
$$
In the special case of ordinary Hopf algebra pairings the mappings
$\mu^l_{A,B}$, etc.\ are given by
$\mu^l_{A,B}(a\otimes b)= b_{(1)}\,\langle a,b_{(2)}\rangle$, etc.
These mappings are actions \cite{Ma2}. The same is true for
multiplier Hopf algebra pre-pairings as it is described in the following
proposition.

\proclaim{Proposition 2.2}
Let $\pair AB$ be a multiplier Hopf algebra pre-pairing.
Then the maps $\mu^l_{A,B}$ and $\mu^r_{A,B}$ are actions of $A$ on $B$,
i.e. $(B,\mu^l_{A,B})$ is an $A$-left module and $(B,\mu^r_{A,B})$
is an $A$-right module respectively. Analogously $\mu^l_{B,A}$
and $\mu^r_{B,A}$ are left and right actions of $B$ on $A$ respectively.
\endproclaim

\demo{Proof}
We use Definition 2.1 and in particular 
$\oml a(\id\otimes\oml {a'})\Delta(b)=\oml {a'}(\oml a\otimes\id)
\Delta(b)=\oml{a\cdot a'}(b)$ for all $a,a'\in A$ and $b, b'\in B$.
Then we arrive at
$$\aligned
\mu^l_{A,B}(a\,a'\otimes b)\cdot b'
&=(\id\otimes\oml {a\,a'})(\Delta(b)(b'\otimes\E))\\
&=(\id\otimes\oml a)\left[(\id\otimes\id\otimes\oml{a'})
  ((\id\otimes\Delta)\Delta(b))\cdot(b'\otimes\E)\right]\\
&=(\id\otimes\oml a)\Delta((\id\otimes\oml{a'})\Delta(b))\cdot b'\\
&=\mu^l_{A,B}(a\otimes\mu^l_{A,B}(a'\otimes b))\cdot b'
\endaligned
\tag 2.5
$$
where we used Lemma 1.4 in the second and third equation. Because
non-degeneracy of the multiplications is assumed eqns.\ (2.5) prove that
$\mu^l_{A,B}$ is an action. In a similar manner all other cases can be
verified.\lfl$\square$
\enddemo
\nl
Henceforth the actions will be denoted by ``$\lact$'' and ``$\ract$''
if the meaning is clear. For example $\mu^l_{A,B}(a\otimes b)=:a\lact b$
and $\mu^r_{B,A}(a\otimes b)=:a\ract b$ which means ``$a$ acts from
the left on $b$" and ``$b$ acts from the right on $a$" respectively,
according to the direction of the arrows ``$\lact$" and ``$\ract$".

\proclaim{Lemma 2.3}
Let $\pair AB$ be a multiplier Hopf algebra pre-pairing. Then we obtain
$$\aligned \langle b\lact a,b'\rangle &=\langle a,b'b\rangle\\
           \langle a\ract b,b'\rangle &=\langle a,b\,b'\rangle\\
           \langle a,a'\lact b\rangle &=\langle a\,a',b\rangle\\
           \langle a,b\ract a'\rangle  &= \langle a'a,b\rangle
\endaligned\tag 2.6
$$
for $a,a'\in A$ and $b,b'\in B$.
\endproclaim

\demo{Proof}
The proof is a direct consequence of Definition 2.1.
For instance we get for $a\in A$ and $b,b'\in B$:
$\langle a,b'b\rangle=\omr {b'b}(a)=\omr{b'}
(\id\otimes\omr b)\Delta(a)=
\langle b\lact a,b'\rangle$.\lfl$\square$
\enddemo
\nl
From Lemma 2.3 we immediately get

\proclaim{Lemma 2.4} If the MHA pre-pairing $\pair AB$ is non-degenerate
then $(A,\mu^l_{B,A},\mu^r_{B,A})$ is a $B$-bimodule and
$(B,\mu^l_{A,B},\mu^r_{A,B})$ is an $A$-bimodule.
\endproclaim

\demo{Proof} Let $a\in A$ and $b_1,b_2,b_3\in B$.
Consider $(b_1\lact a)\ract b_2$ and $b_1\lact (a\ract b_2)$
paired with $b_3$. Using Lemma 2.3 and associativity of $B$
yields ${\langle (b_1\lact a)\ract b_2,b_3\rangle}=
{\langle a,b_2 b_3 b_1\rangle}= {\langle b_1\lact (a\ract b_2),b_3\rangle}$
and this proves the lemma because of the non-degeneracy of the
pre-pairing.\lfl$\square$
\enddemo
\nl
Through the help of Lemma 1.4 the commutation rules of the actions with the
comultiplications of a MHA pre-pairing $\pair AB$ can be determined.
Similarly as in the comments preceding Lemma 1.4 one observes
that for all $a\in A$ and $b\in B$ the following multipliers can be
defined.
$$\left.\aligned
        &(\id\otimes(.)\ract a)\Delta_B(b)\\
        &(\id\otimes a\lact(.))\Delta_B(b)\\
        &((.)\ract a\otimes\id)\Delta_B(b)\\
        &(a\lact(.) \otimes\id)\Delta_B(b) 
        \endaligned
   \right\}
   \in \M(B\otimes B)\,,
   \qquad
   \left.\aligned
        &(\id\otimes(.)\ract b)\Delta_A(a)\\
        &(\id\otimes b\lact(.))\Delta_A(a)\\
        &((.)\ract b\otimes\id)\Delta_A(a)\\
        &(b\lact(.) \otimes\id)\Delta_A(a)
      \endaligned
   \right\}
   \in \M(A\otimes A)\,.
\tag 2.7
$$
For example $[(\id\otimes(.)\ract a)\Delta_B(b)]\cdot
(b'\otimes b''):=
(b_{(1)}b')\otimes (b_{(2)}\ract a)b''$ and
$(b'\otimes b'')\cdot[(\id\otimes(.)\ract a)\Delta_B(b)]
=(b'b_{(1)})\otimes b''(b_{(2)}\ract a)$
for any $b',b''\in B$, where the generalized Sweedler notation is used
which is explained in the Appendix. Hence

\proclaim{Proposition 2.5}
Let $\pair AB$ be a MHA pre-pairing. Then for all 
$a\in A$ and $b\in B$ we have
$$\aligned
  \Delta_B(a\lact b)&=(\id\otimes a\lact(\cdot))\Delta_B(b)\\
  \Delta_B(b\ract a)&=((\cdot)\ract a\otimes\id)\Delta_B(b)
  \endaligned
  \quad
  \text{and}
  \quad
  \aligned
  \Delta^{\roman{op}}_B(a\lact b)&=(a\lact(\cdot) \otimes\id)
  \Delta^{\roman{op}}_B(b)\\
  \Delta^{\roman{op}}_B(b\ract a)&=(\id\otimes(\cdot)\ract a)
  \Delta^{\roman{op}}_B(b)
  \endaligned
\tag 2.8
$$
and analogously for $\Delta_A$.
\endproclaim

\demo{Proof}
Let $a,a'\in A$ and $b,b'\in B$. Using Lemma 1.4 and the
coassociativity of $\Delta$ yields
$$\aligned
\Delta(a\lact b)\cdot(a'\otimes b')&= \Delta((\id\otimes\oml a)\Delta(b))
  (a'\otimes b')\\
&=(\id\otimes\id\otimes\oml a)((\id\otimes\Delta[\Delta(b)(a'\otimes\E)])
  (\E\otimes b')\\
&=(\id\otimes a\lact(\cdot))\Delta(b)\cdot(a'\otimes b')\,.
\endaligned
\tag 2.9
$$
Proceeding in an analogous manner completes the proof of the
proposition.\lfl$\square$
\enddemo

\proclaim{Proposition 2.6}
Let $\pair AB$ be a MHA pre-pairing. Then 
$$\langle T^A_2(a\otimes a'),b\otimes b'\rangle =
         \langle a\otimes a',T^B_1(b\otimes b')\rangle\tag 2.10
$$
for $a, a'\in A$ and $b, b'\in B$. If in addition
$\mu^l_{B,A}$ and $\mu^r_{A,B}$ are surjective, we have
$$\langle S_A(a),b\rangle = \langle a,S_B(b)\rangle\,.\tag 2.11
$$
Analogous results can be derived if $\mu^r_{B,A}$ and $\mu^l_{A,B}$ are
supposed to be surjective. Under the surjectivity condition of the
proposition also the following identities hold.
$$
S^{\pm 1}(b\lact a)=S^{\pm 1}(a)\ract S^{\mp 1}(b)\quad\text{and}\quad
S^{\pm 1}(a\lact b)=S^{\pm 1}(b)\ract S^{\mp 1}(a)\,.\tag 2.12
$$
\endproclaim

\demo{Proof} Let $a,a'\in A$ and $b,b'\in B$. Using Lemma 2.3 and the
generalized Sweedler notation yields
$$\aligned
\langle T_2(a\otimes a'),b\otimes b'\rangle &=\langle (a\otimes\E)
 \Delta(a'), b\otimes b'\rangle\\
 &=\langle a a'_{(1)},b\rangle\langle a'_{(2)},b'\rangle=\langle a\,
 (b'\lact a'),b\rangle\\
 &=\langle b'\lact a', b\ract a\rangle
 \endaligned
 \tag 2.13
$$
and
$$\aligned
\langle a\otimes a',T_1(b\otimes b')\rangle &= \langle a,b_{(1)}\rangle
\langle a',b_{(2)}b'\rangle\\
  &=\langle a',(b\ract a)b'\rangle=\langle b'\lact a',b\ract a\rangle\,.
  \endaligned
  \tag 2.14
$$
which leads to (2.10).
Since $T_1$ and $T_2$ are bijective, eqn.\ (2.10) is also valid
for the inverse mappings. With the help of (1.5) one obtains similarly
as before
$$\align
\langle T_2^{-1}(a\otimes a'),b\otimes b'\rangle &=\langle (S\otimes\id)
  (\Delta(a)(S^{-1}(a')\otimes\E)), b\otimes b'\rangle\\
  &=\langle S((\id\otimes\omr {b'})(\Delta(a)(S^{-1}(a)\otimes\E)),b
  \rangle\\
  &=\langle a\,S(b'\lact a'),b\rangle\\
  &=\langle S(b'\lact a'),b\ract a\rangle\\
  \intertext{and on the other hand}
  &=\langle a\otimes a',T_1^{-1}(b\otimes b')\rangle\\
 &=\langle a',S((\oml a\otimes\id)(\E\otimes S^{-1}(b'))\Delta(b))\rangle\\
  &=\langle a',S(b\ract a)b'\rangle\\
  &=\langle b'\lact a', S(b\ract a)\rangle
\endalign
$$
which proves (2.11) because $\mu^l_{B,A}$ and $\mu^r_{A,B}$ are surjective
by assumption. For the verification of (2.12) we are making use of (2.11).
$$\aligned
 S(b\lact a)\cdot S(a')&=(\id\otimes\omr b)(S\otimes\id)(a'\otimes\E)
 \Delta(a)\\
 &=(\id\otimes\omr{S^{-1}(b)})(\Delta_{\roman{op}}(S(a))(S(a')\otimes\E))\\
 &=(S(a)\ract S^{-1}(b))\cdot S(a')\,.
\endaligned
\tag 2.15
$$
\ \lfl$\square$
\enddemo
\nl
Assuming the surjectivity conditions of Proposition 2.6 we obtain
identities which relate $\mu^r_{A,B}$ and $\mu^l_{A,B}$ in a
multiplier Hopf algebra pre-pairing $\pair AB$.
$$\aligned S_B\left(S_A(a)\lact b\right)S_B(b')
          &=S_B\left(b' (S_A(a)\lact b)\right)\\
          &=S_B\left((\id\otimes \oml{S_A(a)})(b'\otimes\E)\Delta(b)
          \right)\\
          &=(\oml a\otimes\id)\Delta(S_B(b))(\E\otimes S_B(b'))\\
          &=(S_B(b)\ract a)S_B(b')
\endaligned
$$
for any $a\in A$ and $b,b'\in B$.
In a similar manner relations for $\mu^l_{B,A}$ and
$\mu^r_{B,A}$ can be deduced. Explicitely we have

\proclaim{Lemma 2.7} Let $\pair AB$ be a multiplier Hopf algebra
pre-pairing and assume that $\mu^l_{B,A}$ and $\mu^r_{A,B}$
(or $\mu^r_{B,A}$ and $\mu^l_{A,B}$) are surjective. Then it holds
$$\aligned
S^{\pm 1}(b\lact a) &= S^{\pm 1}(a)\ract S^{\mp 1}(b)\\
S^{\pm 1}(a\lact b) &= S^{\pm 1}(b)\ract S^{\mp 1}(a)
\endaligned
\tag 2.16
$$
for any $a\in A$ and $b\in B$.\lfl$\square$
\endproclaim
\nl
With the help of the bracket $\langle\cdot,\cdot\rangle$ of a
multiplier Hopf algebra pre-pairing $\pair AB$ we can define multipliers
according to
$$\gathered R:=(\id\otimes \langle\cdot,\cdot\rangle\otimes\id)\circ
               (\Delta_A\otimes\Delta_B) :A\otimes B\to \M(A\otimes B)\\
            R(a\otimes b)\cdot(a'\otimes b') :=
         (\id\otimes \langle\cdot,\cdot\rangle\otimes\id)
         \left(\Delta(a)(a'\otimes\E)\otimes\Delta(b)
         (\E\otimes b')\right)\\
         (a'\otimes b')\cdot R(a\otimes b) :=
         (\id\otimes \langle\cdot,\cdot\rangle\otimes\id)
         \left((a'\otimes\E)\Delta(a)\otimes(\E\otimes b')
         \Delta(b)\right)
\endgathered\tag 2.17
$$
\nl
where $a,a'\in A$ and $b,b'\in B$. Very analogously we define the mapping
$$\tilde R:=(\id\otimes \langle\cdot,\cdot\rangle\circ\tau\otimes\id)\circ
           (\Delta_B\otimes\Delta_A) :B\otimes A\to \M(B\otimes A)\,.
\tag 2.18
$$

\proclaim{Proposition 2.8}
Let $\pair AB$ be a multiplier Hopf algebra pre-pairing. Then the 
following conditions are equivalent.
\roster
\item $R(A\otimes B)=A\otimes B$.
\item $\mu^l_{B,A}$ is surjective.
\item $\mu^r_{A,B}$ is surjective.
\item $\tilde R(B\otimes A)=B\otimes A$.
\item $\mu^l_{A,B}$ is surjective.
\item $\mu^r_{B,A}$ is surjective.
\endroster
In this case $R:A\otimes B\to A\otimes B$ and
$\tilde R:B\otimes A\to B\otimes A$ are bijective with inverse mappings
$$\aligned
R^{-1}&:=(\id\otimes \langle\cdot,\cdot\rangle\circ(S^{-1}\otimes\id)
  \otimes\id)\circ(\Delta_A\otimes\Delta_B)
   :A\otimes B\to A\otimes B\,,\\
\tilde R^{-1}&:=(\id\otimes \langle\cdot,\cdot\rangle\circ
   (S^{-1}\otimes\id)\circ\tau\otimes\id)\circ(\Delta_B\otimes\Delta_A)
   :B\otimes A\to B\otimes A\,.
\endaligned
$$
\endproclaim

\demo{Proof}
``(1)$\Rightarrow$ (2)'':
Let $a\otimes b\in A\otimes B$. Then by assumption
there is an $\sum_i p_i\otimes q_i\in A\otimes B$ s.t.\
$R(\sum_i p_i\otimes q_i)=a\otimes b$. Hence
applying $(a'\otimes b')\cdot(\cdot)$ and then
$(\id\otimes\oml {a''})$ to both sides of the equation yields
$a'\,a\cdot \oml {a''}(b'\,b)=a'\sum_i r_i\lact p_i$ where
$r_i =(\id\otimes\oml {a''})((\E\otimes b')\Delta(q_i))$. Since $A$
is non-degenerate algebra this yields the result.
\nl
``(2)$\Rightarrow$(1)'': Let $a,a'\in A$ and $b,b',b''\in B$.
Then
$$\align
&(a'\otimes\E)R(b\lact a\otimes b')(\E\otimes b'')\\
&=(\id\otimes\langle\cdot,\cdot\rangle\otimes\id)
  [(a'\otimes\E)\Delta(b\lact a)\otimes\Delta(b')(\E\otimes b'')]\\
&=(\id\otimes\langle\cdot,\cdot\rangle\otimes\id)
  [(\id\otimes\id\otimes\omr b)
 ((\id\otimes\Delta)((a'\otimes\E)\Delta(a))\otimes b'_{(1)}\otimes
 (b'_{(2)}b'')]\\
&=(\id\otimes\omr {(b'_{(1)}b)})((a'\otimes\E)\Delta(a))\otimes
 (b'_{(2)}b'')\\
&=(a'\otimes\E)[(b'_{(1)}b)\lact a\otimes b'_{(2)}](\E\otimes b'')\,.
\endalign
$$
where we used in particular Proposition 2.5 and Definition 2.1. Hence
$$R\circ(\mu^l_{B,A}\otimes\id) = (\mu^l_{B,A}\otimes\id)\circ
   (\tau\otimes\id)\circ(\id\otimes\tau \,T_{\roman{op}\,1}\,\tau)\circ
   (\tau\otimes\id)
\tag 2.19
$$   
from which the statement follows since $\mu^l_{B,A}$ is supposed to be
surjective.
\nl
``(1)$\Leftrightarrow$(3)'': The proof of this equivalence works
pretty similar to the proofs ``(1)$\Leftrightarrow$(2)''.
We consider $a'a\otimes b'b
=(a'\otimes b')R(\sum_j p'_j\otimes q'_j)$ and then we arrive at
the relation
$b'b\,\oml {(a'a)}(b'')={b'\sum_j q'_j\ract((\omr{b''}\otimes\id)
\Delta(p'_j))}$.
On the other side we obtain
from $(a'\otimes\E)R(a\otimes b\ract a'')(\E\otimes b')$
the identity
$$R\circ(\id\otimes\mu^r_{A,B}) = (\id\otimes\mu^r_{A,B})\circ
   (\id\otimes\tau)\circ(\tau\,T_{\roman{op}\,2}\,\tau\otimes\id)\circ
   (\id\otimes\tau)\,.
\tag 2.20
$$   
``(4)$\Leftrightarrow$(5)$\Leftrightarrow$(6)'': 
In an analogous manner the equivalence of the conditions (4), (5) and (6)
can be proved.
\nl
``(3)$\Leftrightarrow$(5)'': Essentially this proof has been done in
(2.16).
\nl
To prove that
$\Theta:=(\id_A\otimes \langle\cdot,\cdot\rangle\circ(S^{-1}\otimes\id)
  \otimes\id_B)\circ(\Delta_A\otimes\Delta_B)
   :A\otimes B\to \M(A\otimes B)$ is the inverse mapping of $R$
observe that
$$\Theta\circ(\mu^l_{B,A,}\otimes\id)=
(\mu^l_{B,A,}\otimes\id)\circ(\tau\otimes\id)\circ
[\id\otimes(S^{-1}\otimes\id)\,T_2\,(S\otimes\id)]\circ(\tau\otimes\id)
\tag 2.21
$$
which can be derived similarly as eqn.\ (2.19). From (1.5) one derives
$$(\tau\otimes\id)\circ(\id\otimes(S^{-1}\otimes\id)\,T_2\,(S\otimes\id))
\circ(\tau\otimes\id)=
[(\tau\otimes\id)\circ(\id\otimes\tau\,T_{\roman{op}\,1}\,\tau)\circ
   (\tau\otimes\id)]^{-1}\,.
$$
Hence the comparison of (2.19) with (2.21) leads to
$\Theta=R^{-1}$ because $\mu^l_{B,A}$ is surjective. Similarly the
inverse of $\tilde R$ can be determined.\lfl$\square$
\enddemo
\nl
We have provided enough results to define the notion of a pairing
of multiplier Hopf algebras.

\definition{Definition 2.9}
A multiplier Hopf algebra pairing $\pair AB$ is a
multiplier Hopf algebra pre-pairing where the conditions of the
previous proposition are fulfilled.
If $A$ and $B$ are multiplier Hopf $*$-algebras then we demand
additionally $\langle a^*,b\rangle=\overline{\langle a,S(b)^*\rangle}$.
\enddefinition

\remark{Remark 3} (1) In particular an ordinary Hopf algebra
pairing $\pair {H_1}{H_2}_{\roman {Hopf}}$ \cite{Ma1,VD1} constitutes
a multiplier Hopf algebra pairing since, for instance, $\mu^l_{H_1,H_2}$
is surjective and therefore the conditions of Proposition 2.8 are
fulfilled. Conversely a non-degenerate MHA pairing $\pair {H_1}{H_2}$ of
the Hopf algebras $H_1$ and $H_2$ is a Hopf algebra pairing. Indeed, 
the (co-)multiplication and antipode properties of
$\langle\cdot,\cdot\rangle$ are obvious from Definition 2.1 and
eq.\ (2.11). From Definition 2.1.(3) we derive $\langle\E,\E\rangle=1$
for $a=a'=\E$ and $b=\E$.
If we put $b=\E$ and $a=\E$ in eq.\ (2.10) and use
the non-degeneracy of $\langle\cdot,\cdot\rangle$ we arrive
at $\langle a'_{(1)},\E\rangle\,a'_{(2)}= a'$.
Applying $\varepsilon$ on both sides of this equation yields
$\langle a',\E\rangle=\varepsilon (a')$. Analogously
$\langle\E, b\rangle=\varepsilon(b)$ is shown.
\nl
(2) If $(A,\Delta)$ is a regular multiplier Hopf ($*$-)algebra with
non-trivial (left) integral $\varphi$, then the dual $\hat A$ of $A$ is
also a regular multiplier Hopf ($*$-)algebra with integral \cite{VD3}.
If we take $\langle\cdot,\cdot\rangle:A\otimes\hat A\to k$ to be the
evaluation map, then $\pair A{\hat A}$ is a non-degenerate MH($*$-)A
pairing. For the verification of this statement we use the results of
\cite{VD3}. Because of \cite{Proposition 3.4, VD3} the algebras
$\hat A=\{\varphi(a\cdot)\mid a\in A\}$ and $A$ are dually paired
vector spaces with respect to the bracket
$\langle a,\varphi(a'\cdot)\rangle:=\varphi(a'a)$. It is obvious that 
$(\omr{\varphi(a'\cdot)}\otimes\id)\Delta(a)=(\varphi\otimes\id)
\left((a'\otimes\E)\Delta(a)\right)\in A$,
and $\omr{\varphi(a\cdot)\cdot\varphi(a'\cdot)} =
\varphi(a\cdot)\cdot\varphi(a'\cdot) 
=\omr{\varphi(a\cdot)}\circ(\id\otimes \omr{\varphi(a'\cdot)})\circ\Delta=
\omr{\varphi(a'\cdot)}\circ(\omr{\varphi(a\cdot)}\otimes\id)\circ\Delta$
by \cite{Proposition 4.2, VD3}, for any $a,a'\in A$.
Similarly the remaining conditions of Definition 2.1 can be proved.
The explicit expression for the action $\mu^l_{\hat A,A}$ is given
through $\varphi(a'\cdot)\lact a= (\id\otimes\varphi)\,
\left((\E\otimes a')\Delta(a)\right)$. Hence the action is surjective
because of the bijectivity of $T_2^{\roman{op}}$. Since
$\varphi(a\cdot)^*(a')=\overline{\varphi(a\cdot S(a')^*)}$ the $*$-property
of the bracket $\langle\cdot,\cdot\rangle$ is a consequence of
\cite{Proposition 4.7, VD3}. Summarizing the results we obtain that
$\pair A{\hat A}$ is a non-degenerate MH($*$-)A pairing.
\endremark

\proclaim{Proposition 2.10} For a non-degenerate multiplier
Hopf algebra pairing $\pair AB$ the actions $\mu^l_{A,B}$, $\mu^r_{A,B}$,
$\mu^l_{B,A}$ and $\mu^r_{B,A}$ are non-degenerate in the sense of
Proposition 1.2.
\endproclaim

\demo{Proof}
For $b\in B$ let $a\lact b=0$ for all $a\in A$. This is equivalent
to $(\id\otimes\oml a)\Delta(b)\cdot b'$ for any $a\in A$ and $b'\in B$.
Acting with $\varepsilon_B$ on both sides and using
$(\id\otimes\varepsilon_B)\,T_{\roman{op}\,1}^B\,(b\otimes b')=
\varepsilon_B(b')\,b$ (see \cite{VD2}) we arrive at
$\varepsilon_B(b')\langle a,b\rangle=0$ for all $a\in A$ and $b'\in B$.   
Since $\varepsilon \ne 0$ it follows $b=0$.
In the same way the non-degeneracy
of all other actions will be proved.\lfl$\square$
\enddemo

\proclaim{Proposition 2.11}
If $\pair AB$ is a multiplier Hopf $*$-algebra pre-pairing,
and for all $a\in A$ and $b\in B$ it holds $\langle a^*,b\rangle=
\overline{\langle a, S(b)^*\rangle}$ and
$\langle a,b^*\rangle=\overline{\langle S(a)^*, b\rangle}$ then 
$$\aligned (a\ract b)^* &= a^*\ract S(b)^*\,,\\
           (b\lact a)^* &= S(b)^*\lact a^*\,,\\
           (b\ract a)^* &= b^*\ract S(a)^*\,,\\
           (a\lact b)^* &= S(a)^*\lact b^*\,.
\endaligned \tag 2.22
$$
\endproclaim

\demo{Proof}
Let $a, a'\in A$ and $b\in B$. Then we obtain
$$\aligned
(a\ract b)^*\cdot a'&=[((\omr b\otimes\id)\Delta(a))^*]\cdot a'
 =((\omr b\otimes\id)(\E\otimes a'{}^*)\Delta(a))^*\\
 &=(\langle a_{(1)},b\rangle a'{}^*a_{(2)})^* = \langle (a^*)_{(1)},S(b)^*
 \rangle a^*_{(2)} a'\\
 &=(a^*\ract S(b)^*)\cdot a'
\endaligned
\tag 2.23
$$
where we used the $*$-property of multipliers according to Chapter 1
in the second equation, the generalized Sweedler notation as explained in
the Appendix in the third equation, and the $*$-property of the
bracket $\langle .,.\rangle$ according to Definition 2.1 in the fourth
equality. Hence the first statement of the proposition is verified.
Similarly all other equations in (2.22) can be proven.\lfl$\square$
\enddemo
\nl
Consider a regular multiplier Hopf ($*$-)algebra $(A,\Delta)$.
It is clear form \cite{VD2,VD3} and the results of the previous chapter
that the opposite co-opposite object
$A^{\roman{op}}_{\roman{op}}:=(A^{\roman{op}},\Delta^{\roman{op}})$
is again a regular multiplier Hopf ($*$-)algebra. $A^{\roman{op}}$
is the opposite algebra to $A$.

\proclaim{Proposition 2.12}
Let $\pair AB$ be a (non-degenerate) multiplier Hopf ($*$-)algebra
pairing. Then $\pair {A^{\roman{op}}_{\roman{op}}}
{B^{\roman{op}}_{\roman{op}}}$
is again a (non-degenerate) multiplier Hopf ($*$-)algebra pairing.
\endproclaim

\demo{Proof}
For $\opop A$ and $\opop B$ we obtain
$$\aligned
[(\id\otimes\omr b)\Delta_{\roman{op}}(a)]_{1^{\roman{op}}}(a')
&=(\omr b\otimes \id)((\E\otimes a')\Delta(a))\\
&=[(\omr b\otimes\id)\Delta(a)]_2(a')
\endaligned
\tag 2.24
$$
and analogous results can be found for the second component.
Because of Definition 2.1 it follows
$$\aligned
\oml a(\id\otimes\oml {a'})\Delta^{\roman{op}}(b) &=
\oml a(\oml {a'}\otimes\id)\Delta(b)=\oml {a'}(\id\otimes\oml a)\Delta(b)
=\oml{a'\,a}\\
&=\oml {a'}(\oml a\otimes\id)\Delta^{\roman{op}}(b)\\
&=\oml {a\diamond a'}
\endaligned
\tag 2.25
$$
where ``$\diamond$'' is the opposite multiplication. If $\pair AB$
is a $*$-pairing then the involution is antimultiplicative w.r.t.\
the opposite multiplication and $\Delta^{\roman{op}}$ is
$*$-homomorphism. The $*$-property of $\langle\cdot,\cdot\rangle$
according to Definition 2.1 holds for $\pair{\opop A}{\opop B}$ since
$\opop S=S$. From relation (2.24) one deduces that the action
$\mu^l_{\opop B,\opop A}$ is surjective because
$\mu^r_{B,A}$ is surjective.\lfl$\square$
\enddemo
\nl
On the tensor product $A^{\roman{op}}_{\roman{op}}\otimes
B^{\roman{op}}_{\roman{op}}$
we can define a mapping according to eqns.\ (2.17).
We will denote it henceforth by $R^{\roman{op}}_{\roman{op}}:=
(\id\otimes \langle\cdot,\cdot\rangle\otimes\id)\circ
(\Delta^{\roman{op}}\otimes\Delta^{\roman{op}})$.
If $\pair AB$ is multiplier Hopf algebra pairing it follows
$R^{\roman{op}}_{\roman{op}}(A\otimes B)=A\otimes B$ (as sets).
From the proof of Proposition 2.8 we obtained particular results
which will be important for further calculations and which we would
like to collect in a lemma.

\proclaim{Lemma 2.13}
Let $\pair AB$ be a multiplier Hopf algebra pre-pairing then the following
identities hold.
$$\aligned
 R(b\lact a\otimes b')&= (b'_{(1)}b)\lact a\otimes b'_{(2)}\\
 R^{-1}(b\lact a\otimes b')&= S^{-1}(S(b)b'_{(1)})\lact a\otimes b'_{(2)}\\
 R(a\otimes b\ract a')&= a_{(1)}\otimes b\ract (a'a_{(2)})\\
 R^{-1}(a\otimes b\ract a')&= a_{(1)}\otimes b\ract S^{-1}(a_{(2)}S(a'))\,.
\endaligned
\tag 2.26
$$
And because of the symmetry reasons outlined in Proposition 2.12
it follows immediately
$$\aligned
 \opop R(a\ract b\otimes b')&= a\ract (b\,b'_{(2)})\otimes b'_{(1)}\\
 {\opop R}^{-1}(a\ract b\otimes b')&= a\ract S^{-1}(b'_{(2)}S(b))
 \otimes b'_{(2)}\\
 \opop R(a\otimes a'\lact b)&= a_{(2)}\otimes (a_{(1)}a')\lact b\\
 {\opop R}^{-1}(a\otimes a'\lact b)&= a_{(2)}\otimes S^{-1}(S(a')a_{(1)})
 \lact b\,.
\endaligned
\tag 2.27
$$
\ \lfl$\square$
\endproclaim
\nl
This lemma will be used for the proof of the next proposition.

\proclaim{Proposition 2.14} Let $\pair AB$ be a non-degenerate multiplier
Hopf algebra pairing. Then the following relations are fulfilled.
$$
\gathered
   R^{\roman{op}}_{\roman{op}}\circ(S^{\pm 1}\otimes S^{\mp 1}) =
      (S^{\pm 1}\otimes S^{\mp 1})\circ R\,,\\
   R\circ R^{\roman{op}}_{\roman{op}}=R^{\roman{op}}_{\roman{op}}
   \circ R\,.
\endgathered\tag 2.28
$$
\endproclaim

\demo{Proof}
For the proof one uses the
first part of Proposition 2.6 and $\Delta^{\roman{op}}\circ S^{\pm 1}
=(S^{\pm 1}\otimes S^{\pm 1})\circ\Delta$ according to the results of
\cite{VD2}. Then
$$\aligned
 \opop R\circ(S^{\pm 1}\otimes S^{\mp 1})(a\otimes b)
  =&(\id\otimes \langle\cdot,\cdot\rangle\otimes\id)((S^{\pm 1}\otimes
  S^{\pm 1})
  \otimes (S^{\mp 1}\otimes S^{\mp 1}))(\Delta\otimes\Delta)(a\otimes b)\\
  =&(S^{\pm 1}\otimes S^{\mp 1})\,R(a\otimes b)\,.
  \endaligned
$$ 
The commutativity of $R$ and $\opop R$ will be proved in two steps.
At first one verifies without problems that
$$\aligned
[R\circ\opop R(a\ract b\otimes b')]\cdot(\E\otimes b'')&=
b'_{(1)}\lact a\ract(b\,b'_{(3)})\otimes (b'_{(2)}b'')\,,\\
[\opop R\circ R(b\lact a\otimes b')]\cdot(\E\otimes b'')&=
(b'_{(1)}b)\lact a\ract b'_{(3)}\otimes (b'_{(2)}b'')\\
\endaligned
\tag 2.29
$$
where $a\in A$ and $b,b',b''\in B$. Now we operate with
$\left(\langle\cdot,(\cdot)c\,c'\rangle\otimes\id\right)$ on both equations
where $c,c'\in B$.
After a little calculation using Lemma 2.3 and keeping the generalized
Sweedler notation in mind, we find
$$\aligned
\left(\langle\cdot,(\cdot)c\,c'\rangle\otimes\id\right)
[R\circ\opop R(a\ract b\otimes b')\cdot(\E\otimes b'')]&=
\langle a\ract b,(b'_{(3)}c)(c'b'_{(1)})\rangle (b'_{(2)}b'')\\
\left(\langle\cdot,(\cdot)c\,c'\rangle\otimes\id\right)
[\opop R\circ R(b\lact a\otimes b')]\cdot(\E\otimes b'')]&=
\langle b\lact a,(b'_{(3)}c)(c'b'_{(1)})\rangle (b'_{(2)}b'')\,.
\endaligned
\tag 2.30
$$
Equations (2.30) prove $R\circ\opop R=\opop R\circ R$ since ``$\lact$"
and ``$\ract$'' are non-degenerate because of
Proposition 2.10.\lfl$\square$
\enddemo

\proclaim{Proposition 2.15}
For a non-degenerate multiplier Hopf $*$-algebra pairing
$\pair AB$ the mapping $R$ and the involution ``$\,*$'' are related
according to
$$\aligned R^{\pm 1}\circ(*\otimes *)&=(*\otimes *)\circ R^{\mp 1}\,,\\
   (R^{\roman{op}}_{\roman{op}})^{\pm 1}\circ(*\otimes *)
   &=(*\otimes *)\circ (R^{\roman{op}}_{\roman{op}})^{\mp 1}\,.
\endaligned
\tag 2.31
$$
\endproclaim

\demo{Proof}
The proof of the proposition is rather straightforward.
Let $a, a'\in A$ and $b,b'\in B$, then
$$\aligned
 (a'\otimes b')\cdot R(a^*\otimes b^*)&=
 (\id\otimes\langle\cdot,\cdot\rangle\otimes\id)
 (\Delta(a)({a'}^*\otimes\E)\otimes\Delta(b)
 (\E\otimes {b'}^*))^*\\
 &=(a_{(1)}{a'}^*)^*\,\overline{\langle a_{(2)}, S^{-1}(b_{(1)})
 \rangle}\otimes
   (b_{(2)}{b'}^*)^*\\
 &=(*\otimes *)[R^{-1}(a\otimes b)\cdot ({a'}^*\otimes {b'}^*)]\\
 &=(a'\otimes b')\cdot (*\otimes *) R^{-1}(a\otimes b)
\endaligned
\tag 2.32
$$
where we used the $*$-property of $\langle\cdot,\cdot\rangle$ and
$\Delta$. The verification of the other cases can be worked out
similarly.\ \ \lfl$\square$
\enddemo
\nl
The two morphisms $R$ and $R^{\roman{op}}_{\roman{op}}$ are the ingredients
for the construction of a twist map which we use in the next chapter
for the definition of the multiplication of the quantum double of a
multiplier Hopf algebra pairing $\pair AB$.
\abs
\abs
\head 3. The Quantum Double\endhead
\abs
\definition{Definition 3.1}
The twist map of a non-degenerate MHA pairing
$\pair AB$ is defined as
$$T:=R\circ(R^{\roman{op}}_{\roman{op}})^{-1}\circ\tau :B\otimes A\to
A\otimes B\,.
\tag 3.1
$$
\enddefinition     
\nl
In the case of Hopf algebra pairings it holds
$T(b\otimes a)= \langle S^{-1}(a_{(1)}),b_{(3)}\rangle\,
 \langle a_{(3)},b_{(1)}\rangle\,a_{(2)}\otimes b_{(2)}$.
This mapping is used in \cite{Ma2,VD1} to construct
the multiplication of the quantum double. And indeed, we will see in the
following that also for an MHA pairing $\pair AB$ the mapping $T$ has
enough properties which enable us to construct a multiplication on
the tensor product $A\otimes B$. Furthermore we can show that even
a multiplier Hopf algebra structure on $A\otimes B$ can be established
which generalizes the quantum double construction of usual Hopf algebra
pairings to the case of MHA pairings.
Before we will prove this fact we have to provide several structural
results. Exploiting Lemma 2.13 we arrive at the following proposition.

\proclaim{Proposition 3.2} Let $\pair AB$ be a non-degenerate
MHA pairing. Then the twist map obeys the relations
$$\gather
 {\align
 T(b''\otimes b\lact a\ract b')&=(b''_{(1)}b)\lact a\ract
 S^{-1}(b''_{(3)}S(b'))\otimes b'_{(2)}\,,\tag 3.2\\
 T(a\lact b\ract a'\otimes a'')&= a''_{(2)}\otimes
 S^{-1}(S(a) a''_{(1)})\lact b\ract (a'a''_{(3)})\,,\tag 3.3\\
 T(a\lact b\otimes b'\lact a')&= (b_{(1)}b')\lact a'_{(2)}\otimes
 S^{-1}(S(a)a'_{(1)})\lact b_{(2)}\,\tag 3.4\\
 T(b\ract a\otimes a'\ract b')&= a'_{(1)}\ract
 S^{-1}(b_{(2)}S(b'))\otimes b_{(1)}\ract(a a'_{(2)})
 \tag 3.5
 \endalign}\\
 {\align
 (a'\cdot(\cdot)\otimes (\cdot)\ract a'')T(b\otimes a)&=
 a'a_{(2)}\otimes S^{-1}(a_{(1)})\lact b\ract(a_{(3)}a'')\,,\tag 3.6\\
 (b'\lact (\cdot)\otimes (\cdot)\cdot b'')T(b\otimes a)&=
 (b'b_{(1)})\lact a\ract S^{-1}(b_{(3)})\otimes b_{(2)}b''
 \tag 3.7
\endalign}
\endgather
$$
for all $a,a',a''\in A$ and $b,b',b''\in B$.
\endproclaim

\demo{Proof}
We use Lemma 2.13 to verify
$$T(b''\otimes b\lact a\ract b')
 =R(b\lact a\ract (b'S^{-1}(b'')_{(1)})\otimes S(S^{-1}(b'')_{(2)}))\,.
\tag 3.8
$$
A short calculation shows that
$$b'S^{-1}(b'')_{(1)}\otimes\Delta(S(S^{-1}(b'')_{(2)}))(b\otimes \E)
  =S^{-1}(b''_{(3)} S(b'))\otimes b''_{(1)} b\otimes b''_{(2)}\,.
\tag 3.9
$$
Inserting (3.9) into (3.8) leads to
$$T(b''\otimes b\lact a\ract b')=({b''}_{(1)}b)\lact a\ract
 S^{-1}(b''_{(3)}S(b'))\otimes b'_{(2)}\,.
$$
which proves (3.2). Analogously identity (3.3) can be verified.
Using Lemma 2.13 and Proposition 2.5 according to
$$\aligned
 T(b\ract a\otimes a'\ract b')
 &=R(a'\ract S^{-1}((b\ract a)_{(2)}S(b'))\otimes (b\ract a)_{(1)})\\
 &=R(a'\ract S^{-1}(b_{(2)}S(b'))\otimes b_{(1)}\ract a)\\
 &=[a'\ract S^{-1}(b_{(2)}S(b'))]_{(1)}\otimes b_{(1)}\ract
  (a[a'\ract S^{-1}(b_{(2)}S(b'))]_{(2)})\\
 &=a'_{(1)}\ract S^{-1}(b_{(2)}S(b'))\otimes b_{(1)}\ract (a a'_{(2)})
\endaligned
$$
yields (3.5). Similar calculations lead to (3.4).
For the proof of (3.6) we consider
$$\aligned
  (a_1\cdot(\cdot)\otimes (\cdot)\ract a_2)\,
  T(a_3\lact b\ract a_4\otimes a)
  &=a_1a_{(2)}\otimes S^{-1}(S(a_3)a_{(1)})\lact(b\ract a_4)\ract
  (a_{(3)}a_2)\\
  &=a_1a_{(2)}\otimes S^{-1}(a_{(1)})\lact(a_3\lact b\ract a_4)\ract
  (a_{(3)}a_2)
\endaligned
$$
where we used (3.3). Since the actions ``$\lact$'' and ``$\ract$''
are surjective we obtain the result. Similarly (3.7) is
shown.\lfl$\square$
\enddemo

\proclaim{Proposition 3.3}
The twist map $T$ and the multiplications $\m_A$ and $\m_B$
obey the following relations.
$$\aligned T\circ(\m_B\otimes\id) &=
           (\id\otimes \m_B)\circ(T\otimes\id)\circ(\id\otimes T)\,,\\
           T\circ(\id\otimes \m_A) &=
           (\m_A\otimes\id)\circ(\id\otimes T)\circ(T\otimes\id)\,.
\endaligned\tag 3.10
$$           
\endproclaim

\demo{Proof}
We prove the first equation. The second one can be derived completely
analogous because of the symmetry of the construction.
From Proposition 3.2 we obtain
$$
T(b'\otimes a\ract b)(\E\otimes b'')={b'}_{(1)}\lact a\ract
 S^{-1}(b'_{(3)}S(b))\otimes (b'_{(2)b'')}
\tag 3.11
$$
and hence
$$\aligned
&[(\id\otimes\m_B)\circ(T\otimes\id)\circ(\id\otimes T)
 (b'\otimes b''\otimes a\ract b)](\E\otimes b''')\\
&=(\id\otimes\m_B)(T\otimes\id)(b'\otimes(b''_{(1)}\lact a)\ract
 S^{-1}(b''_{(3)} S(b))\otimes b''_{(2)}b''')\\
&=(b'b'')_{(1)}\lact a\ract S^{-1}((b'b'')_{(3)}S(b))\otimes
 (b'b'')_{(2)}b'''\\
&= T(b'b''\otimes a\ract b)\cdot(\E\otimes b''')
\endaligned
$$
where we used (3.11) two times.\lfl$\square$
\enddemo
\nl
Thus $T$ behaves like a braiding with respect to
the multiplication and the identity map. Making use of the
properties of $T$ and the associativity of $A$ and $B$, we can therefore
define
an associative algebra on the tensor product
$A\otimes B$ which generalizes the algebra structure of a
quantum double of ordinary Hopf algebras
\cite{Dri,Ma2,VD1} to multiplier Hopf algebras.

\definition{Definition 3.4}
The quantum double $\Cal D\pair AB$ of a non-degenerate
multiplier Hopf algebra pairing $\pair AB$ is the algebra
$(A\otimes B,\m_{\Cal D})$ with the multiplication map defined through
${\m_{\Cal D}:= (\m_A\otimes \m_B)\circ (\id\otimes T\otimes\id)}$.
\enddefinition

\proclaim{Corollary 3.5}
The multiplication $m_{\Cal D}$ of the quantum double is non-degenerate.
\endproclaim

\demo{Proof}
For a fixed $d\in \Cal D$ suppose $d\cdot_{\Cal D}d'=0$ for all
$d'\in\Cal D$. Then $TT^{-1}(d)\cdot_{\Cal D}d'=0$ for all $d'\in \Cal D$.
Because of Proposition 3.3 this is equivalent to
$(\id\otimes\m_B)(T\otimes\id)(\id\otimes\m_A\otimes\id)(T^{-1}(d)\otimes
d')=0$ for any $d'\in \Cal D$. Hence
$T(\id\otimes\m_A)(T^{-1}(d)\otimes a')=0$ for any $a'\in A$
since $\m_B$ is non-degenerate. Thus it follows 
$T^{-1}(d)\cdot(\E\otimes a')=0$ for all $a'\in A$ and therefore $d=0$.
Similarly one proves $d\cdot_{\Cal D}d'=0\ \forall d\in
\Cal D\Leftrightarrow d'=0$.\lfl$\square$
\enddemo

\proclaim{Proposition 3.6}
Let $\pair AB$ be a non-degenerate multiplier Hopf $*$-algebra pairing.
Then $\imath_{\Cal D}:=T\circ(*\otimes *)\circ\tau :A\otimes
B\to A\otimes B$
renders $(\Cal D,\m_{\Cal D},\imath_{\Cal D})$
a non-degenerate $*$-algebra.
\endproclaim

\demo{Proof}
The antilinearity of $\imath_{\Cal D}$ is clear. From Proposition 2.15 we
get
$$\aligned
\imath_{\Cal D}^2 &=T\circ(*\otimes *)\circ\tau\circ T\circ(*\otimes *)
  \circ\tau\\
  &=R\circ{\opop R}^{-1}\circ (*^2\otimes *^2)\circ \opop R\circ R^{-1}\\
  &=\id\,.
\endaligned
\tag 3.12
$$
The antimultiplicativity will be proven as follows. Let $d,d'\in \Cal D$,
then
$$\align
 &(d\cdot_{\Cal D}d')^*\\
 &=\imath_{\Cal D}\circ\m_{\Cal D}(d\otimes d')
  \tag 3.13\\
 &=T\circ\tau\circ(\m_B\otimes\m_A)\circ(\tau\otimes\tau)
 \circ(*\otimes*\otimes*\otimes*)\circ(\id\otimes T\otimes\id)
 (d\otimes d')\\
 &=\m_{\Cal D}\circ(\id\otimes T\,(*\otimes*)\,\tau\,T\otimes
 \id)\circ(\tau\otimes\tau)\circ(\id\otimes\tau\otimes\id)\circ
 (\tau\otimes\tau)\circ(*\otimes\id\otimes\id\otimes*)(d\otimes d')\\
 &=\m_{\Cal D}\circ(\imath_{\Cal D}\otimes\imath_{\Cal D})
 (d\otimes d')\,.\\
 &={d'}^*\cdot_{\Cal D} d^*\,.
\endalign
$$
In the second equation of (3.13) the antimultiplicativity of ``$*$'' is used.
The third identity is derived with the help of Proposition 3.3 and in
the fourth equation we made use of $\imath_{\Cal D}^2=\id$.\lfl$\square$
\enddemo

\remark{Remark 4}
There is no reason why we should prefer $A\otimes B$ instead of
$B\otimes A$ for the construction of the quantum double.
One easily observes that the inverse twist map $T^{-1}:A\otimes B
\to B\otimes A$ obeys analogous relations like (3.8).
If the corresponding quantum double is denoted by $\overline{\Cal D}:=
(B\otimes A, \m_{\overline{\Cal D}})$ it is straightforward to verify
that $T:\overline{\Cal D}\to \Cal D$ is a ($*$-)algebra isomorphism.
\endremark
\abs
We are now investigating how multipliers of $A$ and $B$, and multipliers
of $A\otimes A$ and $B\otimes B$ compose to multipliers of $\Cal D$
and $\Cal D\otimes \Cal D$ respectively. As usual in this chapter we
suppose $\pair AB$ to be non-degenerate multiplier Hopf algebra pairing.

\proclaim{Proposition 3.7}
Let $m\in\M(A)$ , $n\in \M(B)$, $M\in \M(A\otimes A)$ and
$N\in\M(B\otimes B)$ be multipliers. Then $\alpha(m\otimes n)$ defined by
$$\aligned
\alpha(m\otimes n)_1&:=(m_1\otimes\id)\circ T\circ (n_1\otimes\id)
     \circ T^{-1}\\
\alpha(m\otimes n)_2&:=
    (\id\otimes n_2)\circ T\circ (\id\otimes m_2)\circ T^{-1}
    \endaligned
\tag 3.14
$$
is a multiplier in $\M(\Cal D)$, and $\beta(M\otimes N)$ given through
$$\aligned
\beta(M\otimes N)_1&:=(M_1)_{1\,3}\circ (T\otimes T)\circ
                            (N_1)_{1\,3}\circ (T^{-1}\otimes T^{-1})\\
\beta(M\otimes N)_2&:=(N_2)_{2\,4}\circ (T\otimes T)\circ
                            (M_2)_{2\,4}\circ (T^{-1}\otimes T^{-1})
      \endaligned
\tag 3.15
$$
is a multiplier in $\M(\Cal D\otimes\Cal D)$.
$(M_1)_{1\,3},\,(N_1)_{1\,3},\,(M_2)_{2\,4},\,(N_2)_{2\,4}
\in \roman{End}_k(A\otimes B\otimes A\otimes B)$, for instance
$(M_1)_{1\,3}$ operates on the
first and third component as $M_1$.
\endproclaim

\demo{Proof}
We give the proof for $\alpha$. The outlined techniques can be applied
in a similar way for the verification of the statement for $\beta$.
Let $d,d'\in \Cal D$. Then
$$\aligned
&\alpha(m\otimes n)_2(d)\cdot_{\Cal D}d'\\
&=\m_{\Cal D}\circ((\id\otimes n_2)\circ T\circ(\id\otimes\m_2)\circ T^{-1}
 \otimes\id\otimes\id)(d\otimes d')\\
&=(\m_A\otimes\id)(\id\otimes T)(\id\otimes\m_B(n_2\otimes\id)\otimes\id)
  (T\,(\id\otimes m_2)\,T^{-1}\otimes T^{-1})\\
&=(\id\otimes\m_B)(T\otimes\id)(\id\otimes\m_A(m_2\otimes\id)\otimes\id)
  (\id\otimes\id\otimes T\,(n_1\otimes\id))(T^{-1}\otimes T^{-1})
  (d\otimes d')\\
&=d\cdot_{\Cal D}\alpha(m\otimes n)_1(d')
\endaligned
\tag 3.16
$$
In the second and third equation of (3.16) use has been made of
Proposition 3.3. The multiplier property of $m$ and $n$ enters in
the third and fourth equality.
According to the assertions in Chapter 1 this proves the
proposition.\lfl$\square$
\enddemo

\remark{Remark 5}
From Proposition 3.7 it is obvious how to proceed
for higher tensor products. If $M\in \M(A^{\otimes\,n})$ and
$N\in \M(B^{\otimes\,n})$ then for example
$$\gamma(M\otimes N)_2=
 (N_2)_{2,4,\ldots,2n}\circ(\undersetbrace \text{$n$ times}\to
 {T\otimes\ldots\otimes T})\circ
 (M_2)_{2,4,\ldots,2n}\circ(\undersetbrace \text{$n$ times}\to
 {T^{-1}\otimes\ldots\otimes T^{-1}})
$$
is the second component of the multiplier $\gamma(M\otimes N)
\in \M({\Cal D}^{\otimes\,n})$.
\endremark
    
\proclaim{Corollary 3.8} The mappings
$$\aligned
i_A:&\cases \M(A)&\to \M(\Cal D)\\
           m&\mapsto \alpha(m\otimes\E_{\M(B)})
           \endcases
\\
i_B:&\cases \M(B)&\to \M(\Cal D)\\
           n&\mapsto \alpha(\E_{\M(A)}\otimes n)
           \endcases
\endaligned
\quad\text{and}\quad
\aligned
I_A:&\cases \M(A\otimes A)&\to \M(\Cal D\otimes \Cal D)\\
           M&\mapsto \beta(M\otimes\E_{\M(B\otimes B)})
           \endcases
\\
I_B:&\cases \M(B\otimes B)&\to \M(\Cal D\otimes \Cal D)\\
           N&\mapsto \beta(\E_{\M(A\otimes A)}\otimes N)
           \endcases
\endaligned
\tag 3.17
$$
are algebra embeddings. If $\pair AB$ is a non-degenerate
multiplier Hopf $*$-algebra pairing then they are $*$-algebra morphisms.
\endproclaim

\demo{Proof}
We restrict to the proof for $i_A$ because all other cases can be derived
similarly. Looking at the first component of $i_A$ it is obvious that
it is an algebra embedding. Because of the uniqueness of multipliers
coinciding in one of their components it follows that $i_A$ is
an algebra embedding. If $\pair AB$ is MH($*$-)A pairing then we
obtain for any $d\in \Cal D$
$$\aligned
(i_A(m)^*)_1(d)&=(i_A(m)_2(d^*))^*\\
&=T\,\tau\,(*\otimes *)\,T\,(\id\otimes m_2)\,\tau\,(*\otimes *)(d)\\
&=(*\circ m_2\circ *\otimes\id)(d)\\
&=i_A(m^*)_1(d)
\endaligned
\tag 3.18
$$
where in the third equation use has been made of $\imath_{\Cal D}^2=\id$.
Hence $i_A$ is a $*$-algebra morphism\lfl$\square$
\enddemo

\remark{Remark 6} Occasionally we will identify $m$, $n$, $M$ and $N$
with their images under the morphisms of Corollary 3.8. Then it
holds for example $M\cdot N=\beta(M\otimes N)$ for $M\in\M(A\otimes A)$
and $N\in \M(B\otimes B)$.
\endremark

\definition{Definition 3.9}
The comultiplication $\Delta_{\Cal D}$ of the quantum double
$\Cal D$ of a non-degenerate multiplier Hopf algebra pairing $\pair AB$
is defined to be the mapping
$$\Delta_{\Cal D}:=\beta\circ(\Delta_A\otimes\Delta^{\roman{op}}_B)
  : \Cal D\to \M(\Cal D\otimes \Cal D)
 \tag 3.19
$$
where $\beta$ from Proposition 3.7 is used.
\enddefinition

\proclaim{Proposition 3.10}
The linear mappings $T^{\Cal D}_1$, $T^{\Cal D}_2$,
$T^{\Cal D}_{\roman{op}\,1}$ and $T^{\Cal D}_{\roman{op}\,2}:
\Cal D\otimes\Cal D\to \Cal D\otimes\Cal D$
according to Definition 1.3 are bijective.
We obtain the following expressions for
$T^{\Cal D}_1$ and $T^{\Cal D}_2$.
$$\aligned T^{\Cal D}_1&= (T^A_1)_{1\,3}\circ(\id\otimes\id\otimes T)\circ
           (\id\otimes T^B_{\roman{op}\,1}\otimes\id)
           \circ(\id\otimes\id\otimes T^{-1})\,,
           \\
           T^{\Cal D}_2&= (T^B_{\roman{op}\,2})_{2\,4}
           \circ(T\otimes\id\otimes\id)
           \circ (\id\otimes T^A_2\otimes\id)
           \circ(T^{-1}\otimes\id\otimes\id)\,.
\endaligned \tag 3.20
$$
\endproclaim

\demo{Proof}
We outline the proof for $T^{\Cal D}_2$. All other cases can be worked
out in the same fashion. In the proof we use the notation
$$T(b\otimes a) =: \sum_i a^{(i)}\otimes b^{(i)}
\quad\text{and}\quad
T^{-1}(a\otimes b) =: \sum_j b_{(j)}\otimes a_{(j)}\,.
\tag 3.21
$$
Let $a_1,a_2,a_3\in A$ and $b_1,b_2,b_3\in B$. Then
$$\aligned
&T^{\Cal D}_2(a_1\otimes b_1\otimes a_2\otimes b_2)\cdot
(\E_{\Cal D}\otimes a_3\otimes b_3)\\
&=\left[(T\otimes T)\left(\left[(T^{-1}\otimes T^{-1})
  \left(a_1\otimes b_1\otimes a_3\otimes b_3\right)\right]\cdot
  \Delta_A(a_2)_{2\,4}\right)\right]\cdot
  \Delta^{\roman{op}}_B(b_2)_{2\,4}\\
&=\sum T(b_{1\,(i)}\otimes a_{1\,(i)}a_{2\,(1)})(\E\otimes b_{2\,(2)})
   \otimes
  (a_3\otimes b_3)\cdot_{\Cal D} (a_{2\,(2)}\otimes b_{2\,(1)})\,.
\endaligned
  \tag 3.22
$$  
Hence  
$$T^{\Cal D}_2(a_1\otimes b_1\otimes a_2\otimes b_2)
 =\sum T(b_{1\,(i)}\otimes a_{1\,(i)}a_{2\,(1)})(\E\otimes b_{2\,(2)})
  \otimes (a_{2\,(2)}\otimes b_{2\,(1)})
$$
which yields the result.\lfl$\square$
\enddemo

\proclaim{Corollary 3.11}          
$\Delta_{\Cal D}$ is coassociative in the sense of (1.3), i.e.
$$
(T^{\Cal D}_2\otimes\id)\circ(\id\otimes T^{\Cal D}_1)=
(\id\otimes T^{\Cal D}_1)\circ(T^{\Cal D}_2\otimes\id)\,.
\tag 3.23
$$
\endproclaim

\demo{Proof}
Taking the expressions (3.20) for $T^{\Cal D}_1$ and $T^{\Cal D}_2$
and making use of (1.3) for $(A,\Delta_A)$ and $(B,\Delta^{\roman{op}}_B)$
proves the corollary.\lfl$\square$
\enddemo
\nl
Before we will prove in Proposition 3.15 that $\Delta_{\Cal D}$ is a
($*$-)algebra homomorphism, we need three lemmas. We use the notation
(3.21).

\proclaim{Lemma 3.12}
Let $\pair AB$ be a non-degenerate MHA pairing. Then it holds for
$a,a'\in A$ and $b,b'\in B$
$$\align
&\Delta_{\Cal D}(T(b\otimes a))\cdot(\E_{\Cal D}\otimes a'\otimes b')\\
 &=\sum (a^{(i)})_{(1)}\otimes (b^{(i)})_{(2)}\otimes
 ((a^{(i)})_{(2)}\otimes\E)\cdot T((b^{(i)})_{(1)}b'_{(j)}\otimes a'_{(j)})
 \\
\intertext{and} 
 &(\Delta^{\roman{op}}_B(b)\cdot\Delta_A(a))
 \cdot(\E_{\Cal D}\otimes a'\otimes b')\\
 &=\sum (a_{(1)})^{(k)}\otimes (b_{(2)})^{(k)}\otimes
  T(b_{(1)}b'_{(l)}\otimes (a_{(2)}a')_{(l)})\,.
\endalign
$$  
\ \lfl$\square$
\endproclaim

\proclaim{Lemma 3.13} 
Let $a_1,a_2,a_3,a_4\in A$ and $b_1,b_2,b_3,b_4,b_5\in B$.
Then
$$\align
&\left(b_1\lact(\cdot)\otimes b_2\cdot(\cdot)\otimes a_1\cdot(\cdot)
 \otimes\id\right)
 [\Delta_{\Cal D}(T(b_3\ract a_2\otimes a_3\ract b_4))
 \cdot(\E_{\Cal D}\otimes a_4\otimes b_5)]\\
&=b_1\lact a_{3\,(1)}\ract S^{-1}(b_{3\,(3)}S(b_4))\otimes
 b_2 b_{3\,(2)}\otimes a_1 a_{3\,(2)}a_{4\,(2)}\otimes
 [S^{-1}(a_{4\,(1)})\lact b_{3\,(1)}\ract (a_2 a_{3\,(3)}a_{4\,(3)})] b_5
\\
\intertext{and}
&\left(b_1\lact(\cdot)\otimes b_2\cdot(\cdot)\otimes a_1\cdot(\cdot)
 \otimes\id\right) 
 [\Delta^{\roman{op}}_B(b_3\ract a_2)\cdot\Delta_A(a_3\ract b_4)
 \cdot(\E_{\Cal D}\otimes a_4\otimes b_5)]\\
&=(b_1 b_{3\,(2)})\lact a_{3\,(1)}\ract S^{-1}(b_{3\,(4)}S(b_4))\otimes
 b_2 b_{3\,(3)}\otimes a_1 (a_{3\,(2)} a_4)_{(2)}\otimes\\
 &\phantom{=\ \ }\otimes[S^{-1}((a_{3\,(2)}a_4)_{(1)}\lact b_{3\,(1)}\ract
 (a_2(a_{3\,(2)} a_4)_{(3)})] b_5\,.
\endalign  
$$
\ \lfl$\square$
\endproclaim

\proclaim{Lemma 3.14}
Let $a_1,a_2,a_3,a_4\in A$ and $b_1,b_2\in B$. Then the following identity
is fulfilled
$$\aligned
&(b_1 b_{2\,(2)})\lact a_{3\,(1)}\otimes
S^{-1}((a_{3\,(2)}a_4)_{(1)})\lact b_{2\,(1)}\otimes
a_1 (a_{3\,(2)}a_4)_{(2)}\otimes a_2 (a_{3\,(2)}a_4)_{(3)}
\\
&=b_1\lact a_{3\,(1)}\otimes S^{-1}(a_{4\,(1)})\lact b_3\otimes
 a_1 a_{3\,(2)} a_{4\,(2)}\otimes a_2 a_{3\,(3)} a_{4\,(3)}\,.
\endaligned
\tag 3.24
$$
\ \lfl$\square$
\endproclaim
\nl
If one uses Proposition 2.5 and Proposition 3.2 the proofs
of the three lemmas are straightforward, although lengthy calculations
are involved. The first part of Lemma 3.12 has already been shown in
Proposition 3.10. To prove Lemma 3.13 one uses Lemma 3.12. For the
proof of Lemma 3.14 it is convenient to multiply both sides of
(3.24) with some $\left(a_I\cdot(\cdot)\otimes b_I\cdot(\cdot)\otimes
(\cdot)\cdot a_{II}\otimes (\cdot)\cdot a_{III}\right)$ and to verify
this new equality. The non-degeneracy of the multiplication then yields
the identity (3.24). By making use of Lemma 3.13 and Lemma 3.14 we obtain
the important proposition.

\proclaim{Proposition 3.15}
The comultiplication $\Delta_{\Cal D}$ of the quantum double
$\Cal D$ of a non-degenerate multiplier Hopf algebra pairing
$\pair AB$ obeys the identity
$$\Delta_{\Cal D}\circ T(b\otimes a)=\Delta^{\roman{op}}_B(b)
   \cdot\Delta_A(a)\quad\text{for all}\ a\in A,\ b\in B
\tag 3.25
$$
where the identification has been made according to Remark 6.    
Hence $\Delta_{\Cal D}:\Cal D\to \M(\Cal D\otimes\Cal D)$
is an algebra morphism. It is a $*$-algebra
morphism if $\pair AB$ is a multiplier Hopf $*$-algebra pairing.
\endproclaim

\demo{Proof}
Lemma 3.13 and Lemma 3.14 immediately lead to equation (3.25). Then it
is straightforward to prove that $\Delta_{\Cal D}$ is a 
($*$-)algebra homomorphism. We use (3.25), Corollary 3.8 and Remark 6
for the proof.
$$\aligned
\Delta_{\Cal D}((a\otimes b)\cdot_{\Cal D}(a'\otimes b'))
&=\sum \Delta_A(a)\Delta_A({a'}^{(i)})\Delta^{\roman{op}}_B(b^{(i)})
  \Delta^{\roman{op}}_B(b')\\
&=\Delta_A(a)\Delta^{\roman{op}}_B(b)\Delta_A(a')\Delta^{\roman{op}}_B(b')
 \\
&=\Delta_{\Cal D}(a\otimes b)\cdot\Delta_{\Cal D}(a'\otimes b')\,.
\endaligned
\tag 3.26
$$
Similarly one verifies the $*$-property of $\Delta_{\Cal D}$.\lfl$\square$
\enddemo
\nl
Finally we gather the previous results to prove the main theorem
on the construction of a quantum double multiplier Hopf algebra out of
a non-degenerate multiplier Hopf algebra pairing $\pair AB$. This theorem
generalizes the quantum double construction of ordinary Hopf algebra
pairings.

\proclaim{Theorem 3.16}
Let $\pair AB$ be a non-degenerate multiplier Hopf ($*$-)algebra
pairing. Then $(\Cal D, \m_{\Cal D},\Delta_{\Cal D}, (\imath_{\Cal D}))$
is a regular multiplier Hopf ($*$-)algebra. Counit and antipode are
given through $\varepsilon_{\Cal D}=\varepsilon_A\otimes\varepsilon_B$
and $S_{\Cal D}=T\circ\tau\circ(S_A\otimes S^{-1}_B)$ respectively.
If $A$ and $B$ have non-trivial integrals then
$\Cal D$ has a non-trivial integral. Explicitely $\varphi_{\Cal D}
= \varphi_A\otimes\psi_B$ is left integral on $\Cal D$,
if $\varphi_A$ is the left integral on $A$ and $\psi_B$
is the right integral on $B$.
\endproclaim

\demo{Proof}
The previous results show that $(\Cal D,\Delta_{\Cal D},(\imath_{\Cal D}))$
is an MH($*$-)A. From Proposition 3.10 it follows that 
$(\Cal D,\Delta_{\Cal D},(\imath_{\Cal D}))$ is regular. Counit and
antipode of $\Cal D$ are easily determined through the equations \cite{VD2}
$$\aligned
\m_{\Cal D}\circ (T^{\Cal D}_1)^{-1}&=\varepsilon_{\Cal D}\otimes\id\\
\m_{\Cal D}\circ(S_{\Cal D}\otimes\id)&=(\varepsilon_{\Cal D}\otimes
 \id)\circ (T^{\Cal D}_1)^{-1}\,.
\endaligned
$$
If $\varphi_A$ is left integral of $A$ and $\psi_B$ is
right integral of $B$, i.e.\ $(\id\otimes\varphi_A)T^A_2=
\id\otimes\varphi_A$ and ${(\id\otimes\Psi_B)T^B_{\roman{op}\,2}}=
{\id\otimes\Psi_B}$, then it is not difficult to verify
$$(\id\otimes\varphi_A\otimes\psi_B)\circ
  T^{\Cal D}_2= (\id\otimes\varphi_A\otimes\psi_B)\,.
\tag 3.27
$$
Hence $\varphi_A\otimes\psi_B$ is left integral of
$\Cal D$.\lfl$\square$
\enddemo
\abs
\abs
\remark{Acknowledgement} We thank J.\ Kustermans for discussions.
\endremark
\abs\abs
\specialhead{}\centerline{Appendix}\endspecialhead
\abs
We present the ``generalized Sweedler notation''
which is used in the paper.
For a regular multiplier Hopf algebra $(A,\Delta)$ the relation (1.4)
holds. Since
the counit $\varepsilon$ is (non-degenerate) algebra morphism in the sense
of Proposition 1.1 and because of \cite{Theorem 3.6, VD2} one obtains
the identity $(\id\otimes\varepsilon)\circ\Delta=
(\varepsilon\otimes\id)\circ\Delta=\id:\M(A)\to\M(A)$
as in the case of ordinary Hopf algebras. Hence we can
define $\Delta^{(n)}$ for any $n\ge -1$ recursively according to
$$\aligned &\Delta^{(-1)}:=\varepsilon\\
           &\Delta^{(n)} :=(\id\otimes\Delta^{(n-1)})\circ\Delta
                          =(\Delta^{(n-1)}\otimes\id)\circ\Delta
            \quad\text{for all}\ n\ge 0
\endaligned
\tag A.1
$$
using the fact that $\Delta$ is coassociative.
From the definition of $\Delta^{(n)}$ it follows immediately that
$\Delta^{(n)}:\M(A)\to \M(A^{\otimes\,{n+1}})$. We have the following
lemma as a direct consequence of this coassociativity,
resembling the case of ordinary Hopf algebras. 

\proclaim{Lemma A.1} 
$$\Delta^{(n+m+r)} =(\id_{A^{\otimes\,n}}\otimes\Delta^{(m)}\otimes
               \id_{A^{\otimes\,r}})\circ\Delta^{(n+r)}\quad
               \text{for all}\ n,m,r\ge 0\,.
               \tag A.2
$$
\ \lfl$\square$
\endproclaim
\nl
Since the mappings $T_1$, $T_2$, $T_1^{\roman{op}}$ and 
$T_1^{\roman{op}}$ are linear mappings on $A\otimes A$, we obtain

\proclaim{Proposition A.2}
Let $n,m,r\ge 0$ and $a_i\in A$ for $i\in\{1,\ldots,n\}$ and
$a'_j\in A$ for $j\in \{1,\ldots,r\}$. For any $a\in A$,
$\epsilon\in\{-1,1\}$ and $p\in\{1,\ldots, n+m+r\}$ denote
by $a^{(\epsilon,p)}$ the linear mapping
$$a^{(\epsilon,p)}:=
\undersetbrace \text{$p-1$ times}\to{\E\otimes\cdots\otimes\E}\otimes
\lambda_\varepsilon(a)\otimes\!\!\!\!\!\!\!
\undersetbrace \text{$n+m+r-p+1$ times}\to{\E\otimes\cdots\otimes\E}
\in \roman{End}_k(A^{\otimes\,n+m+r+1})
$$ 
where $\lambda_{-1}$ is the left multiplication and $\lambda_1$
is the right multiplication in $A$. Then it holds
$$\aligned
&\Delta^{(n+m+r)}(a)\cdot
a_1^{(\epsilon_1,1)}\cdot\ldots\cdot a_n^{(\epsilon_n,n)}\cdot
a'_1{}^{(\epsilon_{n+1},n+m+1)}\cdot\ldots\cdot a'_r{}^{(\epsilon_{n+r},
n+m+r)}\\
&\in A^{\otimes\,n}\otimes \Delta^{(m)}(A)\otimes A^{\otimes\,r}
\endaligned
\tag A.3
$$
for $\epsilon_1,\ldots,\epsilon_{n+r}\in \{-1,1\}$.
\lfl$\square$
\endproclaim
\nl
This suggests to write symbolically
$\Delta^{(n+m+r)}(a):=a_{(1)}\otimes a_{(2)}\otimes\ldots\otimes
a_{(n+m+r+1)}$.
Then, for instance, we arrive at
$$\align
&\Delta^{(n+m+r)}(a)\cdot
a_1^{(1,1)}\cdot\ldots\cdot a_n^{(1,n)}\cdot 
a'_1{}^{(1,n+m+1)}\cdot\ldots\cdot a'_r{}^{(1,n+m+r)} \tag A.4\\
&=a_{(1)}a_1\otimes\ldots \otimes a_{(n)}a_n \otimes
(a_{(n+1)}\otimes\ldots \otimes a_{(n+m+1)})\otimes
a_{(n+m+2)}a'_1\otimes\ldots \otimes a_{(n+m+r+1)}a'_r
\endalign
$$
and we say that the first $n$ indices and the last $r$ indices are
{\it covered} and the $m+1$ indices in between are {\it free}.
Suppose we restrict to multiplications of the type (A.3)
which guarantee the proper tensor factorization. 
Then we can treat covered indices (in a formal sum) as elements of $A$
and the collection of uncovered indices (in the formal sum) as an element
in $\Delta^{(m)}(A)$. Therefore we can apply
tensor products of morphisms on (A.3) according to this factorization.
These rules are obviously compatible, in particular with the successive
multiplication with another $\Delta^{(n+m+r)}(\tilde a)$ and with another
$\tilde a_1^{(\tilde\epsilon_1,1)}\cdot\ldots\cdot
\tilde a_n^{(\tilde\epsilon_n,n)}\cdot
\tilde a'_1{}^{(\tilde\epsilon_{n+1},n+m+1)}\cdot\ldots\cdot
\tilde a'_r{}^{(\tilde\epsilon_{n+r},n+m+r)}$, because
the multiplier algebra is associative, $\Delta$ is coassociative
algebra morphism, and Lemma A.1 and Proposition A.2 hold.
Analogous results are true for $\Delta^{\roman{op}}$ since it is
also a coassociative algebra morphism.
We call the rules figured out in (A.3) and (A.4) the ``generalized
Sweedler notation'' following the common nomenclature for ordinary
Hopf algebras \cite{Swe}.


\enddocument